\documentclass[11pt]{amsart}

\usepackage{amsmath}
\usepackage{graphicx}
\usepackage{amssymb}
\usepackage{epstopdf}
\usepackage[T1]{fontenc}
\usepackage{tikz}

\newcommand*{\rom}[1]{\expandafter\@slowromancap\romannumeral #1@}
\newcommand{\RNum}[1]{\uppercase\expandafter{\romannumeral #1\relax}}
\newcommand{\Z}{\mathbb{Z}}  
\newcommand{\R}{\mathbb{R}}  

\newcommand{\eps}{\varepsilon}

\newcommand{\F}{\mathbf{F}}

\theoremstyle{plain} 
\newtheorem{thm}{Theorem}[section]
\newtheorem{lem}{Lemma}[section]

\newtheorem{cor}{Corollary}[section]

\newtheorem{qu}{Question}
\newtheorem{prb}{Problem}
\newtheorem{con}{Conjecture}

\theoremstyle{definition}

\theoremstyle{remark}

\title{On product of difference sets for sets of positive density}
\author{Alexander Fish}

\address{School of Mathematics and Statistics, University of Sydney, Australia}
\curraddr{}
\email{alexander.fish@sydney.edu.au}
\thanks{}

\keywords{Difference sets, sum-product estimates}

\subjclass[2010]{Primary: 37A45; Secondary: 11E25, 11T30}

\date{8 February 2017}                                           

\begin{document}
\begin{abstract}
In this paper we prove that given two sets $E_1,E_2 \subset \Z$ of positive density, there exists $k \geq 1$ which is bounded by a number depending only on the densities of $E_1$ and $E_2$ such that $k\Z \subset (E_1-E_1)\cdot(E_2-E_2)$. As a corollary of the main theorem we deduce that if $\alpha,\beta > 0$ then there exist $N_0$ and $d_0$ which depend only on $\alpha$ and $\beta$ 
such that for every $N \geq N_0$ and $E_1,E_2 \subset \Z_N$ with $|E_1| \geq \alpha N, |E_2| \geq \beta N$ there exists $d \leq d_0$ a divisor of $N$ satisfying $d \, \Z_N \subset (E_1-E_1)\cdot(E_2-E_2)$.
\end{abstract}

\maketitle
\section{\textbf{introduction}}
One of the main themes of additive combinatorics is sum-product estimates. It goes back to Erd\"os and Szemer\'edi \cite{ES} who conjectured that for any finite set $A \in \Z$ (or in $\R$), for every $\eps > 0$ we have 
\[
|A+A| + |A \cdot A| \gg |A|^{2-\eps},
\]
where the $A+A = \{a+b \, | \, a,b \in A\}$, and $A\cdot A = \{ a b  \, | \, a, b \in A\}$. Currently the best known estimate is due to Konyagin-Shkredov \cite{KS} and it is based on the beautiful previous breakthrough work by Solymosi \cite{So}: 
\[
|A+A| + |A \cdot A| \gg |A|^{4/3 + c},
\]
for any $c < 5/9813$.

In this paper we study a slightly twisted, but nevertheless related, sum-product phenomenon. Namely, we address the following 

\begin{qu}
For a given \textbf{infinite} set $E \subset \Z$, how much structure does possess the set $(E-E) \cdot (E-E)$?  
\end{qu} 

We will restrict our attention to sets having positive density, see the definition below.

Furstenberg \cite{Fu} noticed a intimate connection between difference sets for sets of positive density, and the sets of return times of a set of positive measure in measure-preserving systems. In this paper we will establish an arithmetic richness of a set of return times of a set of a positive measure to itself within a measure-preserving system. Recall that a triple $(X,\mu,T)$ is a measure-preserving system if $X$ is a compact metric space, $\mu$ is a probability measure on the Borel $\sigma$-algebra of $X$, and $T:X \to X$ is a bi-measurable map which preserves $\mu$.  For a measurable set $A \subset X$ with $\mu(A) > 0$ the set of return times from $A$ to itself is:
\[
R(A) = \{ n \in \Z \, | \, \mu(A \cap T^n A) > 0\}.
\] 
We will denote by $E^2= \{ e^2 \, | e \in E \}$ the set of squares of $E \subset \Z$.
It has been proved by Bj\"orklund and the author \cite{BF} that for any three sets of positive measure $A,B,$ and $C$ in measure-preserving systems there exists $k \ge 1$ (depending on the sets $A,B,$ and $C$) such that $ k \, \Z \subset R(A)\cdot R(B) - R(C)^2$. 
One of the motivations for this work was to show that $k$ in the latter statement depends only on the measures of the sets $A,B,$ and $C$.
We prove the latter, and even more surprisingly, we show that $R(C)$ can be omitted. We have

\begin{thm}\label{main_thm}
Let $(X,\mu,T)$ and $(Y,\nu,S)$ be measure-preserving systems, and let $A \subset X, B \subset Y$ be measurable sets with $\mu(A) > 0,$ and $\nu(B) > 0$. Then there exist $k_0$ depending only on $\mu(A)$ and $\nu(B))$, and $k \le k_0$ such that $k\, \Z \subset R(A) \cdot R(B)$.
\end{thm}

This result has a few combinatorial consequences. To state the first application, we recall that the upper Banach density of a set $E \subset \Z$ is defined by
\[
d^*(E) =\limsup_{N \to \infty} \sup_{a \in \Z} \frac{|E \cap \{a,a+1,\ldots,a+(N-1)\}|}{N}.
\]
Through Furstenberg's correspondence principle \cite{Fu}, we obtain 

\begin{cor}\label{cor1}
Let $E_1,E_2 \subset \Z$ be sets of positive upper Banach density. Then there exist $k_0$ which depends only on the densities  of $E_1$ and $E_2$ and $k \le k_0$ 
such that 
\[
k \, \Z \subset (E_1 - E_1)\cdot (E_2 - E_2).
\]
\end{cor}
\medskip

Another application of Theorem \ref{main_thm} is the following result.

\begin{cor}\label{cor2}
For any $\alpha,\beta > 0$ there exist $N_0$ and $d_0$, depending only on $\alpha$ and $\beta$, such that for every $N \geq N_0$ and $E_1,E_2 \subset \Z_N$ with $|E_1| \geq \alpha N, |E_2| \geq \beta N$ there exists $d \leq d_0$ which is a divisor of $N$ and $d \, \Z_N \subset (E_1-E_1)\cdot(E_2-E_2)$.
\end{cor}
\medskip

Corollary \ref{cor2} implies also that if $p$ is a large enough prime and $E_1,E_2 \subset \Z_p$ satisfy $|E_1| \ge \alpha p, |E_2| \ge \beta p$, then $(E_1 - E_1) \cdot (E_2 - E_2) = \Z_p$. This also follows from a result by Hart-Iosevich-Solymosi \cite{HIS} who proved that if $E \subset \F_{q}$ (where $\F_q$ is a field with $q$ elements) with $|E| \ge q^{3/4 + \eps}$ then for $q$ large enough $(E-E) \cdot (E-E) = \F_q$.

\medskip
\noindent \textit{Acknowledgment:} The work has been carried out during a research visit to Weizmann Institute, Israel. The author would like to thank Feinberg visiting program and Mathematics Department at Weizmann Institute for their support. The author is indebted to Omri Sarig for his constant encouragement and support, Eliran Subag and Igor Shparlinski for enlightening discussions, and Ilya Shkredov for his useful comments on the first version of the paper and for allowing us to reproduce his simplified proof of Lemma \ref{lem1}.
 \medskip
 
\section{\textbf{Proof of Theorem \ref{main_thm}}}

Let us assume that $(X,\mu,T)$ is a measure-preserving system, and let $A \subset X$ be a measurable set with 
$\mu(A) > 0$. Recall that the set of return times of $A$ is defined by
\[
R(A) = \{ n \in \Z \, | \, \mu(A \cap T^n A) > 0 \}.
\]

The theorem will follow from the following  statement.
\medskip


\begin{lem}\label{lem1}
For every $L  \ge 1$ and every $b \in \Z \setminus \{ 0 \}$ there exists $m \le \lfloor \frac{1}{\mu(A)^{L }}\rfloor + 1$ such that
\[
\{mb, 2mb, \ldots, L  m b\} \subset R(A).
\]
\end{lem}

Indeed, let $R(A)$ and $R(B)$ be sets of return times for measurable sets $A$ and $B$ of positive measures. Then choose $N = \lfloor \frac{1}{\nu(B)}\rfloor+1$. Then for every $b \in \Z \setminus \{ 0 \}$ there exist $1 \leq i <  j \le N$ such that $\nu((S^b)^i B \cap (S^b)^j B) > 0$. Then by $S$-invariance of $\nu$ it follows that there exists $1 \leq m \leq N$ ($m = j - i$) such that $mb \in R(B)$.
\medskip

Let us define $L = N!$. By Lemma \ref{lem1}  there exists $n = n(L,\mu(A))$ such that for every $b \in \Z \setminus \{ 0 \}$ there 
exists $m \leq n$ with $\{mb,2mb,\ldots,Lmb\} \in R(A)$. 

Let us define $k = L \cdot n!$. Take any $b \in \Z \setminus \{ 0 \}$. 
By the choice of $n$, there exists $m \le n$ such that $\{mb,2mb,\ldots, Lmb\} \in R(A)$. By the choice of $N$ it follows that there exists $1 \le j \le N$ such that $j \cdot \frac{k}{Lm} \in R(B)$. Also, $\frac{Lm}{j}$ is an integer less or equal than $Lm$, therefore $\frac{Lm}{j} b \in R(A)$. 
Thus $kb = \frac{Lm}{j} b  \cdot j \frac{k}{Lm} \in R(A) \cdot R(B)$. This finishes the proof of Theorem \ref{main_thm}.
\medskip

\noindent \textit{Proof\footnote{This proof of the lemma has been proposed to the author by I. Shkredov.} of Lemma \ref{lem1}.}
Let $(X,\mu,T)$ be a measure-preserving system, and let $A \subset X$ be a measurable set, and let $b \in \Z \setminus \{ 0\}$. We introduce a new product system $Z = \prod_{i=1}^{L } X$ with the transformation $S = \prod_{i=1}^{L } T^{ib}$ and the product measure $\nu = \prod_{i=1}^{L } \mu$. Then $(Z,\nu,S)$ is a measure-preserving system, and the set $\tilde{A} = \prod_{i=1}^{L } A$ has measure
\[
\nu\left( \tilde{A}\right) = \mu(A)^{L } > 0.
\]
Then by Poincar\'e lemma there exists $m \leq \lfloor \frac{1}{\mu(A)^{L }}\rfloor + 1$ such that 
\[
\nu(\tilde{A} \cap S^m \tilde{A}) > 0.
\]
The latter means that for every $1 \le i \le L $ we have
\[
\mu(A \cap T^{ibm} A) > 0.
\]
Therefore, we have $\{bm, 2bm,\ldots,L  b m \} \in R(A)$ for $m \le \lfloor \frac{1}{\mu(A)^{L }}\rfloor + 1$.
%
%

\qed
\medskip

\section{\textbf{Proofs of Corollaries \ref{cor1} and \ref{cor2}}}
 Furstenberg \cite{Fu} in his seminal work on Szemer\'edi's theorem showed:
 \smallskip
 
 \noindent \textbf{Correspondence Principle.} \textit{Given a set $E \subset \Z$ there exists a measure-preserving system $(X,\mu,T)$ and a measurable set $A \subset X$ such that for all $n \in \Z$ we have
 \[
 d^*\left( E \cap (E + n) \right) \ge \mu(A \cap T^n A), 
 \]
 and
 \[
 d^*(E) = \mu(A).
 \]} 
 \medskip
 
\noindent \textit{Proof of Corollary \ref{cor1}.} 
Let $E_1, E_2 \subset \Z$ be sets of positive densities. Then by Furstenberg's correspondence principle there exist measure-preserving systems $(X,\mu,T)$ and $(Y,\nu,S)$ and measurable sets $A \subset X$, $B \subset Y$ that satisfy
\[
\mu(A) = d^*(E_1), \mbox{ } \nu(B) = d^*(E_2)
\]
and 
\[
R(A) \subset E_1 - E_1, \mbox{ } R(B) \subset E_2 - E_2.
\]
By Theorem \ref{main_thm} there exist $k(\mu(A),\nu(B))$ and $k \le k(\mu(A),\nu(B))$ such that $k \, \Z \subset R(A) \cdot R(B)$. The latter statement implies the conclusion of the corollary. 

\qed 

\noindent \textit{Proof of Corollary \ref{cor2}.} Let $\alpha > 0,$ and $\beta > 0$ and let $E_1,E_2 \subset \Z_N$ with $|E_1| \geq \alpha N$, and 
$|E_2| \ge \beta N$. It is clear that $X = \Z_N$ with the shift map $Tx = x + 1(\!\!\!\mod N)$ and the uniform measure $\mu$ on $X$ defined by $\mu(E) = \frac{|E|}{N}$ for any $E \subset X$ is a measure-preserving system. 
 It is also clear that for $(X,\mu,T)$ and the sets $E_1,E_2 \subset X$ we have\footnote{We identify here the ring $\Z_N$ with the set $\{0,1,\ldots,N-1\}$.} $R(E_1) = (E_1 - E_1) + N\, \Z$ and $R(E_2) = (E_2 - E_2) + N \, \Z$. Then by Theorem \ref{main_thm} it follows that if $N \ge N_0$, where $N_0$ depends only on $\alpha$ and $\beta$, then there exist $k(\alpha,\beta)$ and $k \le k(\alpha,\beta)$ such that $k \, \Z \subset R(E_1) \cdot R(E_2)$. Then by the Chinese Remainder theorem for $d = \gcd{(k,N)} \le k$ we have $d \, \Z \subset (E_1-E_1) \cdot (E_2 - E_2) + N \Z$, which implies the statement of the corollary.

\qed  
\medskip

\section{\textbf{Further problems}}

\medskip

To formulate the first problem, we mention a recent result by Bj\"orklund-Bulinski \cite{BB}, who proved, in particular, that for any $E \subset \Z^3$ of positive density there exists $k \ge 1$, depending on the set $E$ and not only on its density, such that 
\[
k \, \Z \subset \{ x^2 - y^2 - z^2 \, | \, (x,y,z) \in E-E\}.
\] 
Recall, the definition of the upper Banach density of a set $E \subset \Z^2$:
\[
d^*(E) = \limsup_{b - a \to \infty, d - c \to \infty} \frac{|E \cap [a,b)\times[c,d)|}{(b - a)(d- c)}.
\]

\begin{prb}
Is it true that given $E_1,E_2 \subset \Z$ of positive density there exist
$k_0$, which depends only on $d^*(E_1)$ and $d^*(E_2)$, and $k \le k_0$ such that 
$k \, \Z \subset (E_1-E_1)^2 - (E_2 - E_2)^2$? If yes, can we show that for any set $E \subset \Z^2$ of positive density there exist $k_0$, which depends only on $d^*(E)$, and $k \le k_0$ such that $k \, \Z \subset \{ x^2 - y^2 \,|\, (x,y) \in E-E\}$?
\end{prb}
\medskip

The next two problems arise naturally by Theorem \ref{main_thm} and the following result proved by Bj\"orklund and the author in \cite{BF}:

\begin{thm}\label{thm3}
Let $E \subset Mat_d^0(\Z) = \{ (a_{ij}) \in \Z^{d \times d} \, | \, tr{(a_{ij})} = 0 \}$ be a set of positive density. Then there exists $k \ge 1$ (which a priori depends on the set $E$ and not only on its density) such that for any matrix $A \in k \cdot Mat_d^0(\Z)$ there exists $B \in E-E$ such that the characteristic polynomial of $B$ coincides with the characteristic polynomial of $A$.
\end{thm}

\begin{prb}
Is it true that given $E \subset \Z^2$ of positive upper Banach density,
there exist $k_0$ that depends only on $d^*(E)$ and $k \le k_0$  such that 
\[
k \, \Z \subset \{ x y \,\,| \,\, (x,y) \in E-E\}?
\]
\end{prb}
\medskip

We also would like to establish the quantitative version of Theorem \ref{thm3}:
\begin{prb} 
Is it true that the parameter $k$ in Theorem \ref{thm3} depends only on the density of the set $E \subset Mat_d^0(\Z)$?
\end{prb}
\medskip

In view of Corollary \ref{cor2} we believe that a similar statement holds true for any finite commutative ring.

\begin{con}
Let $\alpha > 0$. Then there exist $N$ and $k$ depending only on $\alpha$ such that for any finite commutative ring $R$ with $|R| \ge N$ and any set $E \subset R$ satisfying $|E| \ge \alpha |R|$ the set $(E-E)\cdot (E-E)$ contains a subring $R_0$ such that $|R| / |R_0| \le k$. 
\end{con} 

  \end{document}